\documentclass{ajour}

\newcommand{\beq}{\begin{eqnarray}}
\newcommand{\eeq}{\end{eqnarray}}
\newcommand{\bq}{\begin{equation}}
\newcommand{\eq}{\end{equation}}

\newcommand{\Int}{\displaystyle \int}
\newcommand{\Frac}{\displaystyle \frac}

\newcommand{\Sup}{\displaystyle \sup}
\newcommand{\Lim}{\displaystyle \lim}

\newcommand{\R}{-\!\!\!\!\!\!\int}
\def\BbC{{\rm C\!\!\!I}}

\newcommand{\BbR}{{\rm I\!R}}

\def\un { \rm 1\!\!I}

\usepackage{amsmath,amsfonts,amssymb,amsxtra,bbm}
\newtheorem{rem}{Remark}
\newcommand{\be}[1]{\begin{equation}\label{#1}}
\newcommand{\ee}{\end{equation}}

\begin{document}
\thispagestyle{empty}
\large
\parindent 0pt
\parskip6pt
\vspace*{.5in}
%
\journame{Journal of Functional Analysis}
\articlenumber{xxxxxxxx}
\yearofpublication{1999}
\volume{xxxxxx}
\cccline{xxxxxxxx \$25.00}


\authorrunninghead{Dolbeault,  Esteban \&  S{\'e}r{\'e}}
\titlerunninghead{Eigenvalues of operators with gaps.
}


\title{On the eigenvalues of operators with gaps.\\
Application to Dirac operators.}

\author{Jean Dolbeault, Maria J. Esteban and Eric S{\'e}r{\'e}}

\affil{CEREMADE (UMR C.N.R.S. 7534) \\
Universit{\'e} Paris IX-Dauphine \\
Place du Mar{\'e}chal de Lattre de Tassigny \\
75775 Paris Cedex 16 - France }

\email{dolbeaul, esteban or sere@ceremade.dauphine.fr}

\abstract{This paper is devoted to a general min-max
characterization of the eigenvalues in a gap of the essential spectrum
of a self-adjoint unbounded operator. We prove an abstract theorem, then we
apply it to the case of Dirac operators with a Coulomb-like potential. 
The result is optimal for the Coulomb potential.}

\begin{article}
\vskip 0.5cm

\fbox{\parbox{310pt}{
\begin{center}
This paper appeared as\\
{\sc J.~Dolbeault, M.~J. Esteban, and E.~S\'{e}r\'{e}},\\
{\em On the eigenvalues of operators with gaps.\\
{A}pplication to {D}irac operators},\\ J. Funct. Anal.,
  174 (2000), pp.~208--226.\\[4pt]
\textsf{An erratum is appended at the end of the text (2022).}
\end{center}
}}
\vskip 1cm

\noindent {\bf AMS Subject Classification:} 81Q10; 49R05,
47A75, 35Q40, 35Q75

\vskip 0.5cm\noindent {\bf Keywords:} Variational methods, self-adjoint operators,
quadratic forms, spectral gaps, eigenvalues, min-max, Rayleigh-Ritz quotients, Dirac
operators, Hardy's inequality.

\newpage

\section{Introduction}
The aim of this paper is to prove a very general result on the 
variational characterization of the eigenvalues of operators with gaps in the essential
spectrum. 
\medskip
More precisely, let ${\mathcal H}$ be a Hilbert space 
and $A:D(A) \subset {\mathcal H}
 \rightarrow {\mathcal H}$ a
self-adjoint operator. We denote by  ${\mathcal F} (A)$ the form-domain of $A$. 
Let
${\mathcal H}_+$, ${\mathcal H}_-$ be two orthogonal Hilbert subspaces of 
${\mathcal H}$ such that ${\mathcal H}={\mathcal H}_+
{\displaystyle \oplus} {\mathcal H}_-$.
\medskip
We denote $\Lambda_+$, $\Lambda_-$ the projectors on ${\mathcal H}_+$,
 ${\mathcal H}_-$. We assume the
existence of a core $F$ (i.e. a subspace of $D(A)$ which is dense for the
norm $\big
\Vert. \big \Vert_{D(A)}$), such that :
\begin{itemize}
\item[$(i)$] $F_+ = \Lambda_+ F$ and $F_- = \Lambda_- F$ are two subspaces
of ${\mathcal F}(A)$.
\item[$(ii)$] $a=\sup_{x_- \in F_-\setminus \{ 0\}} {(x_-, Ax_-)\over
\Vert x_- \Vert^2_{_{\mathcal H}}} <+\infty $ .\medskip
\end{itemize}

We consider the sequence of min-max levels

\bq \lambda_k = \  \inf_{
 \scriptstyle V \ {\rm subspace \ of \ } F_+  \atop  \scriptstyle {\rm dim}
\ V =
k  } \  \Sup_{  \scriptstyle x \in ( V \oplus F_- ) \setminus \{ 0 \} } \
\Frac{(x, Ax)}{||x||^2_{_{\mathcal H}}} \ ,
\qquad k \geq 1.
\label{min-max} \eq

Our last assumption is
\begin{itemize}
\item[$(iii)$] $\qquad\qquad \lambda_1 > a \ .$
\end{itemize}

Now, let $b = \inf \ (\sigma_{\rm ess} (A) \cap (a, + \infty)) \in [a, +
\infty]$.
For $k\geq 1$, we denote by $\mu_k$ the $k^{\rm th}$ eigenvalue of $A$ in
the interval
$(a, b)$, counted with multiplicity, { \it if this eigenvalue exists}.
If there
is no $k^{\rm th}$ eigenvalue, we take $\mu_k = b$. The main result of this
note is

\begin{theorem}\label{S1}  With the above notations, and under
assumptions $(i)-(ii)-(iii)$,
$$ \lambda_k \ = \ \mu_k \ , \quad \forall\,k \geq 1 \ . $$
As a consequence, $\displaystyle{b\,=\lim_{k\to\infty} \lambda_k\,=\sup_k
\lambda_k\,>\,a\;.}$ 
\end{theorem}

Such a min-max approach was first proposed by Talman \cite{[Ta]} and Datta-Deviah
\cite{[Da-De]} in the particular case  of  Dirac operators with a potential, to
compute numerically their first positive eigenvalue. In that case, the decomposition
of ${\mathcal H}$ was very convenient for practical purposes: each 4-spinor
was decomposed in its upper and lower parts.
Note that in the Physics litterature, other min-max approaches were proposed, for the study
of the eigenvalues of Dirac operators
with a potential  (see for instance \cite{[Dr-Go]}, \cite{[Ku]}).

\vspace{2mm}
A rigorous min-max procedure was then considered by Esteban and S\'er\'e in \cite{[Es-Se]}
for Dirac operators $H_0+V$, $V$ being a Coulomb-like potential. This time, ${\mathcal H}_+$
and ${\mathcal H}_-$ were the positive and negative spectral spaces of the free Dirac 
operator $H_0$.

\vspace{2mm}
To our knowledge, the first abstract theorem on the variational principle 
(\ref{min-max}) is due to Griesemer and Siedentop 
\cite{[Gri-Sie]}. These authors proved an analogue of Theorem \ref{S1}, under
conditions $(i)$, $(ii)$, and two additional hypotheses instead of $(iii)$: they
assumed that
$(Ax, x)>a
\Vert x\Vert^2\,$ for all $x\in
F_+\setminus \{ 0\}$, and they required the operator
$\,(|A|+1)^{1/2}P_-\Lambda_+\,$
to be bounded. Here, $\Lambda_+$ is the orthogonal projection of $\mathcal H$ on 
${\mathcal H}_+$ and $P_-$ is the spectral projection of $A$ for the interval
$(-\infty,a]$,
$i.e.$ $P_-=\chi_{(-\infty,a]}(A)$.

\vspace{2mm}
Then, Griesemer and Siedentop applied their abstract
result to the Dirac operator with potential.
They proved that the  min-max procedure proposed by Talman and Datta-Deviah
was mathematically correct for a particular
class of bounded potentials. In this case, the restrictions on the potentials
were necessary in order to fulfill the requirement $(Ax, x)>a
\Vert x\Vert^2$, $\;\forall\,x\in
F_+\setminus \{ 0\}$. Such a hypothesis excludes the Coulomb
potentials which appear in atomic models. Griesemer and Siedentop
also applied their theorem to the min-max of \cite{[Es-Se]}, but the boundedness of
$\,(|A|+1)^{1/2}P_-\Lambda_+\,$ seems difficult to check in the case of Coulomb
potentials. See the recent work \cite{[Gri-Lew-Sie]}, where this problem
is partially solved.

\vspace{2mm}
In \cite{[Dol-Est-Ser]}, we extended the result of \cite{[Es-Se]} to a larger class of
Coulomb-like potentials and introduced a minimization approach to define the first
positive eigenvalue of $H_0+V$.

\vspace{2mm}
The present work is motivated by the abstract result of Griesemer and Siedentop
\cite{[Gri-Sie]}. Our Theorem \ref{S1} contains, as particular cases, the results
on the min-max principle for the Dirac
operator of \cite{[Es-Se]}, \cite{[Gri-Sie]}, \cite{[Dol-Est-Ser]}, \cite{[Gri-Lew-Sie]}.
It also applies to the Talman and Datta-Deviah
procedure for atomic Coulomb potentials, under optimal conditions.
However, Griesemer-Siedentop's abstract result is not a consequence of
Theorem \ref{S1}.
Indeed, their hypothesis $(Ax_+, x_+)>a
\Vert x_+\Vert^2\quad(\forall\,x_+\in
F_+\setminus \{ 0\})\;$ does not imply $(iii)$.

\vspace{2mm}
In Section 2 of this paper we prove Theorem \ref{S1}. The arguments are based on
an abstract version of those in \cite{[Dol-Est-Ser]} (\S 4: the minimization
procedure).

\vspace{2mm}
When appling Theorem \ref{S1} in practical situations, the main difficulty
is to check assumption $(iii)$. For that purpose,
an abstract continuation principle (Theorem \ref{TT5})
will be given in Section 3.

\vspace{2mm}
In Section 4 we use Theorems \ref{S1} and \ref{TT5} 
to justify two variational procedures for the eigenvalues of Dirac operators
$H_0+V$: first, Talman's and Datta-Deviah's procedure; then, the min-max
principle of \cite{[Dol-Est-Ser]}. In both cases we cover a large
class of potentials $V$ including Coulomb potentials
$-Z\alpha/|x|$, as long as $Z\alpha<1$. This
condition is optimal since it is well-known that when $Z\alpha\to 1^-$, the
first eigenfunction ``disappears".
For each min-max, we obtain new Hardy-type inhomogeneous inequalities as
by-products of the proof.

\section{ Proof of Theorem \ref{S1}.}

The inequality $\lambda_k\leq \mu_k$ is an easy consequence of conditions $(i)$ and $(ii)$
(see \cite{[Gri-Sie]} for the proof in a similar situation).
It remains to
prove that $\lambda_k\geq \mu_k$ for all $k$. The additional assumption $(iii)$
will be needed, but for the moment, we only assume $(i)$ and $(ii)$.

\vspace{1mm}
 We recall
the notation $\displaystyle{a=\sup_{x_- \in F_-\setminus \{ 0\}} {(x_-,
Ax_-)\over
\Vert x_- \Vert^2_{_{\mathcal H}}} <+\infty} $.
For $E > a$ and $x_+ \in F_+$, let us define
\begin{eqnarray*}
\varphi_{E, x_+} : && F_- \rightarrow \BbR \\
&& y_- \mapsto \varphi_{E, x_+} (y_-) = \Bigl( (x_+ + y_-), A(x_+ + y_-)
\Bigr) - E
||x_+ + y_- ||^2_{_{\mathcal H}} . \end{eqnarray*}
>From assumption $(ii)$, $N (y_-) = \sqrt{(a+1)||y_-||^2_{_{\mathcal H}} -
(y_-, Ay_-)}$ is a norm on $F_-$. Let
$\overline{F}_-^{^N}$ be the completion of $F_-$ for this norm. Since
$||.||_{_{\mathcal H}}
\leq N$ on $F_-$, we have $\overline{F}_-^{^N} \subset {\mathcal H}_-$. For all $x_+\in F_+$,
there is an $x\in F$ such that
$\Lambda_+ x = x_+\,$. If we consider the new variable $z_-=y_--\Lambda_-x$, we can
define 
$$\psi_{E, x}(z_-):= \varphi_{E, \Lambda_+x}(z_-+\Lambda_-x)=(A(x+z_-), x+z_-)-E(x+z_-, x+z_-)\,.$$
Since $F$ is a subspace of $D(A)$, $\psi_{E, x}$ (hence $\varphi_{E, x_+}$) is
well-defined
and continuous for $N$, uniformly on bounded sets. So, $\varphi_{E, x_+}$ has a unique
continuous extension
$\overline \varphi_{E, x_+}$ on $\overline{F}_-^{^N}$, which is continuous for
the extended norm $\overline N$. It is well-known (see e.g. \cite{[Ree-Sim2]}) that
there is a unique self-adjoint operator
$B : D(B) \subset {\mathcal H}_- \rightarrow {\mathcal H}_-$ such that $D(B)$ is a subspace of
$\overline{F}_-^{^N}$, and
\bq \overline {N} (x_-)^2 \ = 
\ (a+1)|| x_-||^2_{_{\mathcal H}} + (x_-, B x_-) \ ,
\quad \forall\,x_- \in D(B) \ . \label{defB}\eq
Now, $\overline \varphi_{E, x_+}$ is of class $C^2$ on
$\overline{F}_-^{^N}$ and
\begin{eqnarray}\label{coer} D^2 \overline\varphi_{E, x_+} (x_-) \cdot (y_-,y_-) & = & -
2 (y_-, By_-) - 2E  ||
y_- ||_{\mathcal H}^2\nonumber \\
& \leq & - 2 \min \left(1, (E-a) \right)\, \overline N (y_-)^2 \ . 
\end{eqnarray}
So $\overline \varphi_{E, x_+}$ has a unique maximum, at
the point $y_- = L_E (x_+)$. The Euler-Lagrange equations
associated to this maximization problem are :
\bq\label{N7} \Lambda_-Ax_+-(B+E)y_-=0\,.\eq

In the sequel of this note,
we shall use the notation ${\mathcal X}'$ for the dual of
a Hilbert space ${\mathcal X}$. Note that $(B+E)^{-1}$ is well-defined and
bounded from ($\overline{F}_-^{^N})^{^{'}}$ to $\overline{F}_-^{^N}$,  since
$E>a$ and $(y_-, (B+a) y_-) \geq 0$, $\;\forall\,y_- \in D(B)$. Moreover,
$x_+ \in F_+ =
\Lambda_+ F$ is of the form $x_+ = \Lambda_+ x = x-\Lambda_- x$ for some $x
\in F \subset
D(A)$. By assumption $(i)$, $\Lambda_- x \in {\mathcal F}(A) \cap
(\overline{F}_-^{^N})$,
and finally $\Lambda_-A(x_+) =\Lambda_- Ax-\Lambda_- A\Lambda_-x\in
(\overline{F}_-^{^N})^{^{'}}$. So the expression $(B+E)^{-1}\Lambda_-x_+$ is
meaningful, and we have 
\bq\label{a1} L_E = (B+E)^{-1} \ \Lambda_- A \ .\eq

\begin{rem}\label{N8} 
{\it The unique maximizer of $\overline\psi_{E, x}:=\overline \varphi_{E,\Lambda_+x}
(\cdot+\Lambda_-x)$ is the vector
$z_-= M_Ex:=L_E\Lambda_+x-\Lambda_-x$ and one has the
following equation for $M_Ex$:
\bq\label{N9} M_Ex=(B+E)^{-1}\Lambda_-(A-E)x\,.\eq
This expression is well-defined, since $x\in D(A)$.}
\end{rem}
\vspace{2mm}

The above arguments allow us, for any $E > a$, to define a map
\begin{eqnarray}
\label{N5} Q_E : & F_+ \rightarrow & \BbR \nonumber\\
& x_+  \mapsto & Q_E (x_+)  = \sup_{x_- \in F_-} \ \varphi_{_{E, x_+}} (x_-) =
\overline{\varphi}_{E, x_+} (L_E x_+) \\
&&\hspace{3mm} =\left(x_+, (A-E)x_+\right)+\left(\Lambda_-Ax_+,(B+E)^{-1}\Lambda_-Ax_+\right)
\;.\nonumber
\end{eqnarray}
Note that for any $x\in F$,
\begin{eqnarray}\label{N5bis} & Q_E (\Lambda_+x)=(x, Ax) & + \ 2\,\hbox{Re }(Ax, M_E x)
\nonumber\\
&& - \ (M_E x, BM_E x)-
E\Vert x+M_E x\Vert ^2
\;.\end{eqnarray}

\vspace{1mm}

It is easy to see that $Q_E$ is a quadratic form with domain $F_+ \subset
{\mathcal H}_+$.

We may also, for $E > a$ given, define  the norm

\bq\label{ne} n_E (x_+) =
||x_+ + L_E x_+ ||_{_{\mathcal H}}\; . \eq 

The following lemma gives some useful inequalities
involving $n_E$ and $Q_E$, and a new formulation of $(iii)$ :

\begin{lemma} \label{N1} Assume that $(i)$ and $(ii)$ are satisfied.  If $\,a<E<E'$,
then
\bq\label{N2bis} \Vert\cdot\Vert_{\mathcal H}\leq n_{E'}\leq n_E\leq
\frac{E'-a}{E-a}n_{E'}\;,\eq 
\bq\label{N2}(E'-E)n_{E'}^2\leq Q_{E}-Q_{E'}\leq (E'-E)n_E^2\;.\eq
Moreover, for any $\,E>a\,$ : 
$$\lambda_1>E\quad\hbox{if and only if}\quad Q_E(x_+)>0\,,\quad\forall\, x_+\in
F_+\;.$$
$$\lambda_1\geq E\quad\hbox{if and only if}\quad Q_E(x_+)\geq 0\,,\quad\forall\,
x_+\in F_+\;.$$

As a consequence, $(iii)$ is equivalent to
\begin{itemize}
\item[$(iii')$] $\quad$ For some $\;E>a\, ,\quad$ $Q_E(x_+)\geq 0\,,\quad\forall\,\,x_+\in F_+\,.$
\end{itemize}
\end{lemma}
\noindent{\it Proof.} Inequality (\ref{N2bis}) is easily proved using the spectral
decomposition of $B$, the formula
$$n_E(x_+)^2=||x_+||^2_{_{\mathcal H}}+||(B+E)^{-1}\Lambda_-Ax_+||^2_{_{\mathcal H}}$$
and the standard inequality
$$1\leq {t+u\over t+v}\leq {u\over v}\;,\quad
\forall\,t\geq 0\,,\; u\geq v >0\;.$$
On the other hand, (\ref{N2}) is a consequence of
$$Q_{E'}(x_+)\geq {\overline\varphi}_{E', x_+}(L_E(x_+))\,, \hbox{ for all }\;E,
E'>a\;.$$
Finally, the definition of $\lambda_1$
implies that
$Q_{E}(x_+)>0$ for all $x_+ \in F_+\setminus\{0\}$ and $a<E<\lambda_1$. But
(\ref{N2bis}) and (\ref{N2})
imply that
$$Q_{\lambda_1}(x_+)\geq Q_{E}(x_+)+
(E-\lambda_1){(\lambda_1-a)^2\over (E-a)^2}\, n^2_{\lambda_1}
(x_+)\;.$$
Passing to the limit $E\to\lambda_1$, we obtain $Q_{\lambda_1}(x_+)\geq 0\;.$

\vspace{1mm}
In the case $\,E>\lambda_1\,$, it follows from the definition of
$\lambda_1$ that for some $\,x_+\in F_+\setminus\{0\}\,$ and some $\varepsilon>0$,  
$$(x_++x_-, A(x_++x_-))\leq (E-\varepsilon)||x_++x_-||^2\;,\quad\forall\,\,x_-\in F_-\;.$$
Hence
$$\varphi_{E, x_+}(x_-)\leq -\varepsilon ||x_++x_-||^2\;,\quad\forall\,\,x_-\in F_-$$
and $Q_E(x_+)\leq -\varepsilon ||x_+||^2<0\;.$ This ends the proof of Lemma
\ref{N1}.

\hfill $\square$

We are now going to give a new definition of the numbers $\lambda_k$, equivalent
to formula
(\ref{min-max}). First of all, let us recall the standard definitions and results on
Rayleigh-Ritz quotients (see e.g. \cite{[Ree-Sim4]}).

Let $T$ be a self-adjoint operator on a Hilbert space $X$, with domain
$D(T)$ and
form-domain ${\mathcal F}(T)$. If $T$ is bounded from below, we may define a
sequence of min-max levels,
$$ \ell_k (T)  \ = \  \inf_{{
 \scriptstyle Y \ {\rm subspace \ of \ } {\mathcal F}(T)  \atop  \scriptstyle
{\rm dim} \
Y = k  }} \  \Sup_{{  \scriptstyle x \in Y\setminus \lbrace 0 \rbrace }}
\ \Frac{(x,
Tx)}{||x||^2_{_{X}}}  . $$
To each $k$ we also associate the (possibly infinite) multiplicity number
$$m_k (T)  \ = \ {\rm card} \Bigl\{ k'
\geq 1 \ , \ \ell_{_{k'}} (T) = \ell_k (T) \Bigr\} \ \geq \ 1 \ .$$
Then $\ell_k (T) \leq \inf  \sigma_{\rm ess} (T)$.
In the case $\ell_k (T) < \inf  \sigma_{\rm ess} (T)$, $\ell_k$ is an
eigenvalue of $T$ with multiplicity $m_k (T)$.
\vspace{2mm}

As a consequence, if ${\mathcal C} \subset {\mathcal F} (T)$ is a form-core for $T$ (i.e.
a dense
subspace of ${\mathcal F} (T)$ for $||.||_{_{{\mathcal F} (T)}}$), then there is a
sequence $(Z_n)$
of subspaces of ${\mathcal C}$, with dim $(Z_n) = m_k (T)$ and
$$ \Sup_{{  \scriptstyle z \in Z_n \atop \scriptstyle ||z||_{_X} = 1 }}  \
|| Tz -
\ell_k (T) z ||_{_{({\mathcal F} (T))'}} \quad \longrightarrow_{_{_{ \hspace{-6mm} n
\rightarrow \infty \hspace{4mm}}}} 0 \ . $$

\vspace{1mm}

Coming back to our situation, we consider the completion
$X$ of $F_+$ for the norm $n_E$.
By (\ref{N2bis}), $X$ does not depend on $E>0$. We denote by $\overline n _E$
the extended norm, and by $<\cdot ,\cdot>_E$ its polar form:
$$<x_+,x_+>_E=(\overline n_E(x_+))^2\,,\;\forall
x_+\in X\,.$$
Since $n_E (x_+)
\geq ||x_+||_{_{\mathcal H}}$, $X$ is a subspace of ${\mathcal H}_+$.  

\vspace{4mm}
{\it We now assume that $(iii)$ is satisfied, i.e. $\,\lambda_1>a$.}
We may define another norm on $\,F_+$ by
$${\mathcal N}_E (x_+) = \sqrt{Q_E (x_+) + (K_E+1)
(n_E (x_+))^2}\, $$
with $K_E = \max \left(0, {(E-a)^2 (E-\lambda_1)\over(\lambda_1-a)^2}\right)\;.$

>From (\ref{N2bis}) and (\ref{N2}), 
${\mathcal N}_E$ is well-defined and  satisfies ${\mathcal N}_E
\geq n_E$. Indeed, 
in the case $a < E \leq \lambda_1$, Lemma \ref{N1} implies $Q_E (x_+) \geq
0$ for all $x_+ \in F_+$.
When $E \geq \lambda_1$,  again from Lemma \ref{N1}, we have
\bq\label{N2ter}
Q_E  \geq Q_{\lambda_1}+(\lambda_1-E)\,n_{\lambda_1}^2\geq -K_E\,n_E^2\,.\eq

Note that for any $a<E<E'$, Lemma \ref{N1} implies the existence of two positive
 constants,
$\,0<c(E, E')<1< C(E, E')\,$, such that
\bq\label{equiv}c(E, E')\,{\mathcal N}_{E'}
\leq {\mathcal N}_E\leq C(E, E')\,{\mathcal N}_{E'}\,.\eq

Let us consider the
completion $\,G\,$ of $\,F_+\,$ for the norm $\,{\mathcal N}_{E}$. Since 
${\mathcal N}_E\geq n_E$, $G$ is a subspace of $X$, dense for the extended norm $\bar
n_E$.
>From (\ref{equiv}),
$G$ does not depend on $E$.
The
extension $\bar Q_E$ of $Q_E$ to $G$ is a closed quadratic form with
form-domain
$G$. So (see e.g. \cite{[Ree-Sim2]}) there is a unique self-adjoint operator $T_E \ : \
D(T_E)\subset X \to X$
with form-domain ${\mathcal F} (T_E) = G$, such that $\bar Q_E (x_+) = <x_+,
T_E x_+>_{E}$, for
any $x_+ \in D(T_E)$. Then $F_+$ is a form-core of $T_E$. The min-max levels $\ell_k(T_E)$
are given by
\bq\label{rayl} \ell_k (T_E)  \ = \  \inf_{{
 \scriptstyle V \ {\rm subspace \ of \ } G  \atop  \scriptstyle
{\rm dim} \
Y = k  }} \  \Sup_{{  \scriptstyle x_+ \in V\setminus \lbrace 0 \rbrace }}
\ \Frac{\bar Q_E (x_+)}{(\bar n_E(x_+))^2}\;  . \eq

The next lemma explains the relashionship between $\ell_k(T_E)$ and the min-max principle
(\ref{min-max}) for $A$.

\begin{lemma} \label{N0} Under assumptions $(i)$, $(ii)$, $(iii)$ :

{\bf (a)} for any $x_+\in F_+\setminus \{0\}$, the real number
$$\lambda(x_+):=\Sup_{  \scriptstyle x \in ( Span(x_+) \oplus F_- ) \setminus \{ 0 \} } \
\Frac{(x, Ax)}{||x||^2_{_{\mathcal H}}}$$
is the unique solution in $\,(a, +\infty)\,$ of the
nonlinear equation
\bq\label{fix1} Q_\lambda (x_+)=0\;.\eq
This equation may be written
\bq\label{fix2} \lambda \Vert x_+\Vert^2_{_{\mathcal H}}=(x_+,Ax_+)+(\Lambda_-Ax_+,(B+\lambda)^{-1}
\Lambda_-Ax_+)\; .\eq

{\bf (b)} The min-max principle (\ref{min-max}) is equivalent to
\bq \lambda_k = \  \inf_{
 \scriptstyle V \ {\rm subspace \ of \ } F_+  \atop  \scriptstyle {\rm dim}
\ V =
k  } \  \Sup_{  \scriptstyle x_+ \in  V \setminus \{ 0 \} } \
\lambda(x_+) \ ,
\qquad k \geq 1 .
\label{min-max2} \eq

{\bf (c)} For any $k\geq 1$, the level $\lambda_k$ defined by (\ref{min-max}) is the
unique solution in $\,(a, +\infty)\,$ of the nonlinear equation
\bq \ell_k  ( T_{\lambda}) \ = \ 0 \ . \label{2} \eq
In other words, $0$ is the $k^{th}$ min-max level for the Rayleigh-Ritz quotients of
$T_{\lambda_k}$, and this determines $\lambda_k$ in a unique way. Moreover, for $a<\lambda\neq
\lambda_k$, the signs of $\lambda_k-\lambda$ and $\ell_k ( T_{\lambda} )$ are the same.
\end{lemma}

\newpage\noindent{\it Proof.} 

{\bf (a)} From Lemma \ref{N1}, $Q_\lambda(x_+)$ is a decreasing continuous function of $\lambda$,
such that $Q_{\lambda_1}(x_+)\geq 0\,$ and $\Lim_{\scriptstyle\lambda\to +\infty}
Q_\lambda(x_+)=-\infty\;.$
So the equation $Q_\lambda(x_+)=0$ has one and only one solution $\tilde{\lambda}(x_+)$, which lies in
the interval $[\lambda_1,+\infty)\;.$ Equation (\ref{fix2}) is equivalent to (\ref{fix1})
by easy calculations. Now, if $\lambda<\tilde{\lambda}(x_+)$, then $Q_\lambda(x_+)>0$,
hence $\lambda(x_+):=\Sup_{  \scriptstyle x \in ( Span(x_+) \oplus F_- ) \setminus \{ 0 \} } \
\Frac{(x, Ax)}{||x||^2_{_{\mathcal H}}} >\lambda\;.$ Similarly, $\lambda>\tilde{\lambda}(x_+)\,$
implies $\lambda (x_+) <\lambda\;.$ So we get
$$\tilde{\lambda}(x_+)=\lambda(x_+)\;.$$

{\bf (b)} Since $\lambda(x_+)=\!\!\!\!\!\!\Sup_{  \scriptstyle x \in  Span(x_+\!) \oplus F_-  \atop
x\ne 0  } \!\! \Frac{(x, Ax)}{||x||^2_{_{\mathcal H}}}\;,$ (\ref{min-max}) is obviously equivalent
to (\ref{min-max2}).

{\bf (c)} We follow the same arguments as in the proof of (a). From Lemma \ref{N1},
the map $\lambda\to \ell _k(T_\lambda)$ is continuous, and
$\ell _k(T_{\lambda_1})\geq 0\;,$ $\Lim_{\scriptstyle\lambda\to +\infty} \ell_k(T_{\lambda})=-\infty\;.$
As a consequence, the equation $\ell_k(T_\lambda)=0$
has at least one solution $\tilde{\lambda}_k$ which lies
in the interval $[\lambda_1,+\infty)\;.$ Now, if $\lambda<\tilde{\lambda}_k\,$ then
from Lemma \ref{N1},
$\ell_k(T_\lambda) >0\;$.
Hence $\Sup_{  \scriptstyle x \in ( V \oplus F_- ) \setminus \{ 0 \} } \
\Frac{(x, Ax)}{||x||^2_{_{\mathcal H}}} >\lambda\,$ for any k-dimensional subspace $V$ of
$F_+$. Similarly, $\lambda>\tilde{\lambda}_k\,$
implies $\Sup_{  \scriptstyle x \in ( V \oplus F_- ) \setminus \{ 0 \} } \
\Frac{(x, Ax)}{||x||^2_{_{\mathcal H}}} <\lambda\,$ for some k-dimensional subspace $V$ of
$F_+$. So, we get
$\tilde{\lambda}_k=\lambda_k\;.$

\hfill $\square$

\vspace{2mm}

As already mentioned, $F_+$ is a form-core of $T_E$ and $G$ is its form-domain.
>From Lemma \ref{N0} (c), $\lambda_k =
\lambda_{k'}$ if and only if $l_{k'}(T_{\lambda_k})=0$. So, denoting
$m_k := {\rm card} \ \lbrace k' \geq 1 \ ; \ \lambda_k =
\lambda_{k'} \rbrace$,  there is a sequence $(Z_n)$ of subspaces of
$F_+$, of dimension $m_k$, such that
$$  \Sup_{{  \scriptstyle x_+ \in Z_n \atop \scriptstyle ||x_+ ||^2_{_{\mathcal H}} +
|| L_{\lambda_k}
x_+ ||^2_{_{\mathcal H}} = 1}} \  \big \Vert T_{\lambda_k} x_+ \big
\Vert_{_{G'}}
\quad \longrightarrow_{_{_{ \hspace{-7mm} n \rightarrow
\infty \hspace{4mm}}}} 0  \; .$$
Using the explicit expressions of $Q_E$ and $L_E$ on $F_+$ (see (\ref{a1}),(\ref{N5bis})),
we obtain
\bq\label{N10}\quad\; \Sup_{{  \scriptstyle \tilde x \in (\un + L_{\lambda_k}) (Z_n)
\atop
\scriptstyle
|| \tilde x ||_{_{\mathcal H}} = 1}} \  \Sup_{\tilde y \in (\un + L_{\lambda_k}) (F_+)
\atop \tilde y\ne 0} \
\Frac{\left|{\mathcal A}(\tilde x, \tilde y)-\lambda_k(\tilde x, \tilde y)_{\mathcal H}
\right|}
{\left((K_{\lambda_k}+1) \ || \tilde y
||^2_{_{\mathcal H}} + Q_{\lambda_k} (\Lambda_+ \tilde y)\right)^{1/2}}   
\longrightarrow_{_{_{ \hspace{-6mm} n \rightarrow  \infty \hspace{1mm}}}} 0
\ , \eq
where $\,{\mathcal A}(\tilde x, \tilde y):= (x, Ay)+(Ax, M_{\lambda_k}y)+(M_{\lambda_k}x,
Ay)-(BM_{\lambda_k}x, M_{\lambda_k}y)$, with $x, y\in F\subset D(A)$ such that $
\Lambda_+x=\Lambda_+
\tilde x$,
$ \Lambda_+y=\Lambda_+ \tilde y$ and\break
$\,M_{\lambda_k}x=L_{\lambda_k}\Lambda_+x-\Lambda_-x\,$. Note that the value of 
$\,{\mathcal A}(\tilde x, \tilde y)\,$ does not depend on the choice of $x$ and $y$. Indeed,
$\,{\mathcal A}\,$ is the polar form of the quadratic form $\,\tilde y\mapsto
{\overline Q}_{\lambda_k}(\Lambda_+\tilde y)+\lambda_k||\tilde y||^2_{\mathcal H}$.

\vspace{1mm}
Denote $\tilde Z_n = (\un + L_{\lambda_k}) (Z_n)$. Take $y \in F$, and let
$\tilde y = (\un + L_{\lambda_k}) (\Lambda_+ y)$. There is a constant
$C(\lambda_k)$
such that
\bq\label{N6}
 (K_{\lambda_k}+1) \ || \tilde y
||^2_{_{\mathcal H}} + Q_{\lambda_k} (\Lambda_+ y) 
\leq  C(\lambda_k) \ || y ||^2_{_{D(A)}} \ . \eq   
Indeed, by Remark \ref{N8}, 
\begin{eqnarray}
&Q_{\lambda_k}(\Lambda_+y) = \left( (A-\lambda_k)y,y+M_{\lambda_k}y\right) \hspace{6cm}\nonumber\\
 &\leq (1\!+\!|\lambda_k|) ||y||_{_{D(A)}}\left(||y||_{_{\mathcal
H}}\!+\!||M_{\lambda_k}y||_{_{\mathcal H}}\right) \! \leq\!(1+|\lambda_k|) \left(1+
{1\!+\!|\lambda_k|\over\lambda_k-a}\right) ||y||^2_{_{D(A)}}.\nonumber
\end{eqnarray}

Moreover, for any $ x \in F_+$, and any $z_-
\in {\mathcal F}(B)$, by (\ref{N9})
we have : \newline
$\left ( (Ax-BM_{\lambda_k}x) - \lambda_k (x +M_{\lambda_k}x), 
z_-\right ) = 0$.
As a consequence, (\ref{N10}) is equivalent to
$$\Sup_{{  \scriptstyle \tilde x \in  \tilde Z_n \atop
\scriptstyle ||
\tilde x
||_{_{\mathcal H}} = 1}} \  \sup_{y \in F\setminus \{ 0\}} \ \Frac{\left|
\left( \tilde x, Ay-\lambda_k y\right)\right|}
{||y||_{_{D(A)}}}  \quad \longrightarrow_{_{_{ \hspace{-6mm} n
\rightarrow
 \infty \hspace{4mm}}}} 0 \ . $$

So, by the standard spectral theory of self-adjoint operators, we obtain an
alternative: either $\lambda_k \in
\sigma_{{\rm ess}} (A) \cap (a, +\infty)$, or $\lambda_k$ is an
eigenvalue of $A$ in the interval $(a,+\infty)$,
with multiplicity greater than or equal to $m_k$.

We have thus proved the inequality $\lambda_k \geq \mu_k$,
$\;\forall\,k \geq 1$. This ends the proof of Theorem \ref{S1}.
\hfill $\square$

\section{An abstract continuation principle.}

This section is devoted to a general method for checking condition $(iii)$ of
Theorem \ref{S1}. It applies to $1$-parameter families of self-adjoint operators
of the form $\,A_\nu=A_0+{\mathcal V}_\nu\,$, with ${\mathcal V}_\nu\,$
bounded. The idea is to prove $(iii)$ for all
$\,A_\nu\,$ knowing that one of them
satisfies it, and having spectral information on every $\,A_\nu$.

\vspace{2mm}
More precisely, we start with a self-adjoint operator  $A_0:D(A_0) \subset {\mathcal H}
 \rightarrow {\mathcal H}\;.$ We denote by  ${\mathcal F} (A_0)$ the form-domain of
$A_0$. 

\vspace{1mm}
For $I$ an interval containing $0\,,$ let $\,\nu\mapsto {\mathcal V}_\nu\,$ a  map
whose values are bounded self-adjoint operators and which is continuous for the
usual norm of bounded operators
$$\displaystyle{\,|||{\mathcal V}|||=\sup_{x\in
{\mathcal H}\setminus\{0\}}\frac{||{\mathcal V} x||_{\mathcal H}}{||x||_{\mathcal H}}\;.}$$
In order to have consistent notations, we also assume that ${\mathcal V}_0=0$. 

\vspace{1mm}
Since $A_0$ is self-adjoint and ${\mathcal V}_\nu$ symmetric and bounded, the operator
$\,A_\nu\,$ is self-adjoint with $\,{\mathcal D}(A_\nu)=
{\mathcal D}(A_0)$,   $\,{\mathcal F}(A_\nu)=
{\mathcal F}(A_0)$.
Let
${\mathcal H}={\mathcal H}_+\oplus {\mathcal H}_-$ be an orthogonal splitting
of ${\mathcal H}$, and $\Lambda_+\;,\;\Lambda_-$ the associated projectors, as
in Section 1.
We  assume the existence of a core
$F$ (i.e. a subspace of $D(A_0)$ which is dense for the norm $\big
\Vert. \big \Vert_{D(A_0)}$), such that :
\begin{itemize}
\item[$(j)$] $F_+ = \Lambda_+ F$ and $F_- = \Lambda_- F$ are two subspaces
of ${\mathcal F}(A_0)$.\medskip

\item[$(jj)$] There is $\,a_-\in\BbR$ such that for all $\,\nu\in
I$, $$a_\nu:={\sup_{x_-
\in F_-\setminus
\{ 0\}} {(x_-, A_\nu x_-)\over
\Vert x_- \Vert^2_{_{\mathcal H}}} \leq a_-}\,. $$
\end{itemize}

For $\nu \in I$, let $b_\nu:=\inf (\sigma_{\rm ess}(A_\nu)\cap (a_\nu,+\infty))\,,$
and for $k\geq 1$, let $\mu_{k,\nu}$ be the $k$-th eigenvalue of $A_\nu$
in the interval $(a_\nu,b_\nu)$, counted with multiplicity, if it exists.
If it does not exist, take $\mu_{k,\nu}:=b_\nu\;.$ Our next assumption is

\begin{itemize}
\item[$(jjj)$] There is $\,a_+>a_-\,$ such that for all 
$\,\nu\in
I$, $\mu_{1,\nu}\geq a_+\,.$
\end{itemize}

Finally, we define the levels
\bq \lambda_{k,\nu} := \  \inf_{
 \scriptstyle V \ {\rm subspace \ of \ } F_+  \atop  \scriptstyle {\rm dim}
\ V =
k  } \  \Sup_{  \scriptstyle x \in ( V \oplus F_- ) \setminus \{ 0 \} } \
\Frac{(x, A_\nu x)}{||x||^2_{_{\mathcal H}}} \ ,
\qquad k \geq 1\, ,
\label{min-max-nu} \eq
and our last assumption is

\begin{itemize}
\item[$(jv)$]  $\quad \lambda_{1,0}>a_-\;$.
\end{itemize}

The main result of this section is

\begin{theorem}\label{TT5} Under conditions $(j)$ to $(jv)$, $A_\nu$ satisfies the
assumptions $(i)$ to $(iii)$ of Theorem \ref{S1} for all $\,\nu\in I$, and 
$\,\lambda_{k,\nu}= \mu_{k,\nu}\geq a_+$, for all $\,k\geq 1$.
\end{theorem}
Note that the boundedness assumption on $\,{\mathcal V}_\nu\,$ is rather restrictive.
However, as it will be seen in Section 4, unbounded perturbations can also be dealt
with, thanks to a regularization argument.

\vspace{2mm}
\noindent {\it Proof of Theorem \ref{TT5}}. Assumptions $(i)$, $(ii)$ of Theorem
\ref{S1} are of course satisfied for all $\,\nu\in I$ : see $(j)$, $(jj)$. From
formula (\ref{min-max-nu}), it is clear that for all $\nu,\nu'\in I\,,$
$$ |\lambda_{1,\nu}-\lambda_{1,\nu'}|\,\leq\,|||{\mathcal V}_{\nu}-
{\mathcal V}_{\nu'}|||\;.$$
So the map $\;\nu\in I\to \lambda_{1,\nu}\;$ is continuous. The set
$$P:=\{\nu\in I : \lambda_{1,\nu}\geq a_+\}$$
is thus closed in $I$, and the set
$$P':=\{\nu\in I : \lambda_{1,\nu}> a_-\}$$
is open. Obviously, $P\subset P'\,.$ But if $\nu\in P'$ then $A_\nu$ satisfies
$(iii)$, so it follows from Theorem \ref{S1} that
$$\,\lambda_{k,\nu}= \mu_{k,\nu}\geq a_+\;, \quad \hbox{for all}\;\,k\geq 1\;,$$
hence $\nu\in P\,.$ As a consequence, $P=P'$, and $P$ is open and closed in $I\,.$
But $P$ is nonempty : it contains $0$. So, $P$ coincides with $I$.
\hfill $\square$

\section{Applications and remarks :  Dirac operators.}

With the notations of the preceding sections, let us define $\,{\mathcal H} = L^2
(\BbR^3,
\BbC^{\,4}),$

\vspace{1mm}
Let $F=C^\infty_0(\BbR^3,\BbC^{\,4})$ be the space of smooth, compactly
supported functions from $\BbR^3$ to 
$\BbC^{\,4}.$

The free Dirac operator is ${ H}_0 = -i \alpha \cdot \nabla +
\beta\,,$ with
$$ \alpha_1, \alpha_2, \alpha_3, \beta \in {\mathcal M}_{4 \times 4} (\BbC\,), \
\beta =
\left ( \begin{array}{cc} \un & 0 \\
0 & - \un \\ \end{array} \right )  \ , \ \alpha_i= \left (
\begin{array}{cc} 0 &
\sigma_i \\ \sigma_i & 0 \\ \end{array} \right ) \ , $$
$\sigma_i$ being the Pauli matrices
$$\sigma _1=\begin{pmatrix}0 & 1\\ 1 & 0 \end{pmatrix}  \quad ,\quad  \sigma_2=\begin{pmatrix} 0 & -i \\
i & 0 \end{pmatrix} \quad ,\quad  \sigma_3=\begin{pmatrix} 1 & 0\\ 0 &-1\end{pmatrix} \, .$$
Let $V$ be a scalar potential satisfying 
\bq\label{V1} V(x) \quad\;\longrightarrow_{_{{\hspace{-8mm}}|x|\to +\infty}} \; 0 \,, \eq
\bq\label{V2} -\frac{\nu}{|x|}-c_1\leq V\leq c_2 =\Sup (V)\,,
 \eq 
with $\,\nu\in (0,1)$, $c_1,c_2\in\BbR$.

\vspace{1mm}
Under the above assumptions, $ H_0 + V$ has a distinguished self-adjoint extension $A$
with domain ${\mathcal D}(A)$ such that
$$ H^1 (\BbR^3, \BbC^4) \subset {\mathcal D} (A) \subset 
H^{1/2} (\BbR^3, \BbC^4) \,, $$
$$\sigma_{\rm ess
} (A) \ = \ (- \infty, - 1] \cup [1, +
\infty) \ , $$
and $F$ is a core for $A$ (see \cite{[Th]}, \cite{[Schmin]},\cite{[Ne]}, \cite{[Kl-Wu]} ). 
In the sequel,
we shall denote this extension indifferently by $A$ or $H_0+V$.
We shall also denote $\,\mu_k(V)\,$
the k-th eigenvalue of $H_0+V$ in the interval $(c_2-1, 1)$, with the understanding
that $\mu_k(V)=1$ whenever $H_0+V$ has less than $k$ eigenvalues in $(c_2-1, 1)$.

\vspace{2mm}

In this section, we shall prove the validity of two different variational
characterizations of the eigenvalues $\mu_k(V)$  corresponding to two different
choices of the splitting $\,{\mathcal H}={\mathcal H}_+\oplus{\mathcal H}_-$,
under conditions which are optimal for the Coulomb potential. In both cases, this
will be done using Theorem \ref{S1}. The main difficulty is to check assumption $(iii)$
of this theorem. It will be sufficient to do it for the Coulomb potential
$V_\nu:=-\nu/|x|\,$. Then, by a simple comparison argument, all potentials satisfying
(\ref{V1}), (\ref{V2})  with the additional condition

\bq\label{V3}\quad c_1, c_2\geq 0,\;\, c_1+c_2-1< \sqrt{1-\nu^2}
\eq 
will be covered by our results. The constant $\sqrt{1-\nu^2}$ is the smallest
eigenvalue of $H_0-\frac{\nu}{|x|}$ in the interval $(-1,1)$.
\vspace{1mm}

The Coulomb potential is not bounded. In order to apply Theorem \ref{TT5}, we shall
use a regularization argument.
 The method will be the following: first replace $\,V_\nu=-\frac{\nu}{|x|}\,$ by
$\,V_{\nu, \varepsilon}:=-\frac{\nu}{|x|+\varepsilon}\,$, $\varepsilon>0$. Then apply
Theorem
\ref{TT5} to
$\,A_{\nu, \varepsilon}:=H_0+V_{\nu, \varepsilon}\,$, for $\varepsilon>0$
fixed and $\nu$ varying in $I=[0,1)$, and $a_+=0,\,a_-=-1$. Combined with Lemma
\ref{N1}, this theorem gives
$$\,Q_{0,\nu,
\varepsilon}(x_+)\geq 0\;,\quad \forall\,\,x_+\in F_+$$
where, following (\ref{N5}), 
\begin{eqnarray}\label{N5-nu-eps}
&\hspace{-1cm}\displaystyle{Q_{E,\nu,
\varepsilon}(x_+)\!:= \!\!\! \sup_{x_-\in F_-}\!\Bigl( (x_+ + y_-), A_{\nu,
\varepsilon}(x_+ + y_-)
\Bigr)\! - E
||x_+ + y_- ||^2_{_{\mathcal H}} }\nonumber
 \\
&
 =\left(x_+, (A_{\nu,
\varepsilon}-E)x_+\right)+\left(\Lambda_-A_{\nu,
\varepsilon}x_+,(B_{\nu, \varepsilon}+E)^{-1}\Lambda_-A_{\nu,
\varepsilon}x_+\right)\;,
\nonumber
\end{eqnarray}
and 
$B_{\nu,\varepsilon}:D(B_{\nu,\varepsilon})\subset {\mathcal H_-}\to {\mathcal H_-}$ is a self-adjoint operator
such that $(x_-,A_{\nu,\varepsilon} x_-)=-(x_-,B_{\nu,\varepsilon} x_-)$ for all $x_-\in F_-\,$ : see \S 2, formula
(\ref{defB}).\vspace{2mm}
Passing to the limit $\,\varepsilon\to 0\,$ in the above inequality, we get
$$Q_{0,\nu, 0}(x_+)\geq 0\,, \quad\forall\,x_+\in F_+,$$
and by Lemma \ref{N1}, this is equivalent to assumption $(iii)$ of Theorem \ref{S1} for
the operator $H_0-{\nu\over |x|}\,.$

\vspace{2mm}
\subsection{The min-max of Talman and Datta-Deviah.}

In this subsection, we choose the following splitting of $\,{\mathcal H}$ :

 $$
{{\mathcal H}^T_+ = L^2
(\BbR^3,
\BbC^{\,2}) \otimes \left\{\left(0\atop 0\right)\right\}}\,, \;{{\mathcal H}^T_- =\left\{\left
(0\atop 0\right)\right\} \otimes L^2
(\BbR^3,
\BbC^{\,2})}\,,$$
so that, for any $\psi=\left(\varphi\atop\chi\right)\in
L^2(\BbR^3,\BbC^{\,4})$, $$\Lambda^T_+\psi=\left(\varphi\atop 0\right),\, 
\Lambda^T_-\psi=\left(0\atop\chi
\right).$$

 With this choice, let $\lambda^T_k(V)$ be the $k$-th min-max associated to $\,A=H_0+V$
by formula (\ref{min-max}). In the case $k=1$, we have

\bq\label{Talm}\lambda^T_1(V)=\inf_{\varphi\neq 0}\sup_{\chi}{(\psi,(H_0+V)\psi)\over
(\psi,\psi)}\;.\eq

This is exactly the min-max principle of Talman (\cite{[Ta]}) and
Datta-Deviah (\cite {[Da-De]}).
It is clear that under conditions (\ref{V1})- (\ref{V2}), assumptions $(i)$ and $(ii)$ of
Theorem
\ref{S1} are satisfied, with
$$a=\sup_{x_- \in F_-\setminus \{ 0\}} {(x_-, Ax_-)\over
\Vert x_- \Vert^2_{_{\mathcal H}}}=c_2-1\;.$$

\vspace{1mm}
The main result of this subsection is

\begin{theorem}\label{TT1} Let $V$ a scalar potential satisfying
(\ref{V1})-(\ref{V2})-(\ref{V3}).
Then,  for all $\,k\geq 1$,  
\bq\label{18}\lambda^T_k( V)=\mu_k( V)\,.\eq
Moreover, 
 $\lambda^T_k(V)= \mu_k(V)\,$ is given by
\bq\label{L11} \lambda^T_k (V)= 
 \inf_{{Y  \ {\rm subspace \ of \ } C^\infty_{o} (\BbR^3,
\BbC^{\,2})\atop{ \atop{\rm dim } Y=k}}} \ \sup_{\varphi\in
Y\setminus\{0\}}\lambda^{T}\!(V,\varphi)\,,\eq where

\bq \lambda^{T}\!(V,\varphi):=\sup_{\psi=\left(\varphi\atop\chi\right)\atop \chi\in
{C^\infty_0} (\BbR^3, \BbC^{\,2})}\frac{((H_0+V)\psi, \psi)}{(\psi,\psi)}\label{XX1}\eq
is the unique number in $(c_2-1,+\infty)$ such that 
\bq \lambda^T(V,\varphi) \Int_{\BbR^3} |\varphi|^2 dx \! = \! \Int_{\BbR^3} \Bigl(
\Frac{|
(\sigma\cdot\nabla) \varphi|^2}{1-V+\lambda^T(V,\varphi)} + (1+V) |\varphi|^2
\Bigr) dx
\label{LA} \eq
The maximizer of (\ref{XX1}) in ${\mathcal H}^T_-$ is
\bq\label{XXX1} \chi(V,\varphi):= \frac{-i(\sigma\cdot\nabla)\varphi}
{1-V+\lambda^T(V,\varphi)}\,.\eq

\end{theorem}

\begin{rem}\label{mini}
{\it In the case $k=1$,  the min-max (\ref{L11}) reduces to
$$\lambda^T_1 (V)= \inf_{\varphi\in C^\infty_{o} (\BbR^3,
\BbC^{\,2})\setminus\{0\}} \lambda^{T}\!(V,\varphi)\;,$$
where $\lambda^{T}(V,\varphi)$ is given by equation (\ref{LA}).
This formulation is equivalent to the minimization principle of \cite{[Dol-Est-Ser]},
\S 4, formula (4.16).}
\end{rem}

\vspace{1mm}
\noindent {\it Proof of Theorem \ref{TT1}}. 

Formulas (\ref{L11}), (\ref{LA}), (\ref{XXX1}) are simply those of Lemma \ref{N0}
(a)-(b), rewritten in the context of the present subsection. So the only thing to prove
is (\ref{18}). For that purpose, we just have to check that condition $(iii)$ of
Theorem \ref{S1} is fulfilled by $H_0+V$. In view of Remark \ref{mini}, this was already
done in \cite{[Dol-Est-Ser]}. But the arguments can be made simpler and clearer, thanks
to the formalism of Sections 2 and 3.
\vspace{1mm}

First of all, since $\,\lambda_1\,$ is monotonic in $V$, it is
sufficient to check $(iii)$ when $\,V_\nu=-\frac{\nu}{|x|}\,$, for all $\nu\in [0, 1)$.

The key inequality that we use below is the following :

\bq\label{import} \mu_1(V)\geq 0\quad \hbox{as soon as}\quad 
-\frac{\nu}{|x|}\leq V\leq 0\;,\;0\leq \nu <1\,.\eq
This inequality can be found in \cite{[wust]}. In the particular case
of Coulomb potentials, it is well-known that
\bq
\,\mu_1(-\frac{\nu}{|x|})=
\sqrt{1-\nu^2}\quad\mbox{for}\quad 0\leq\nu<1\;.
\label{wust2}
\eq

We proceed in two steps.

\vspace{2mm}
{\sl First step :} for $\,\nu\in I:=[0, 1)\,$ and $\varepsilon\geq 0\,$, let
$\, V_{\nu, \varepsilon}:=-\frac{\nu}{|x|+\varepsilon}\;$.
We now fix $\varepsilon>0\,$. The one-parameter family $\,\nu\in I\to
A_{\nu, \varepsilon}:= H_0+V_{\nu, \varepsilon}\,$ and the projectors $\Lambda^T_\pm$
satisfy all the assumptions of Theorem \ref{TT5},
with $\,a_-=-1\,$ and $\,a_+=0$. In particular, $(jjj)$ follows from
(\ref{import}). So we obtain
$$\lambda^T_1(V_{\nu, \varepsilon})=\mu_1(V_{\nu, \varepsilon})\geq
0\,,$$
for all $\,\nu\in [0, 1)$. From Lemma \ref{N1}, this can be written as
\bq\label{1et} Q^{T}_{0,\nu, \varepsilon}(\varphi)\geq 0\,,\quad \forall\,\varphi\in
C^\infty_0(\BbR^3, \BbC^2)\,,\eq
with 

\bq\label{QVet} 
Q^{T}_{E,\nu, \varepsilon}(\varphi) \ = \ \Int_{\BbR^3} \Bigl(
\Frac{|
(\sigma\cdot\nabla) \varphi|^2}{1+E-V_{\nu, \varepsilon}} + (1-E+V_{\nu, \varepsilon}) |\varphi|^2
\Bigr) dx \;.\eq

\vspace{2mm}

{\sl Second step :} For $\,\nu\in [0, 1)\,$ and $\,\varphi\in C^\infty_0(\BbR^3, \BbC^2)
\,$ fixed, we pass to the limit $\,\varepsilon\to 0\,$ in (\ref{1et}). We get :
\bq\label{2et}
Q^T_{0,\nu,0}(\varphi) \geq 0\,,\quad \forall\,\varphi\in
C^\infty_0(\BbR^3, \BbC^2)\,.\eq

So $A_{\nu,0}=H_0+V_\nu$ satisfies criterion $(iii')$ of Lemma \ref{N1}, which is
equivalent to $(iii)$. By Theorem \ref {S1}, we thus have
$$\lambda_1^T(V_\nu)= \mu_1(V_\nu)=\sqrt{1-\nu^2}\,,$$
for all $\,\nu\in (0, 1)$. This ends the proof.
\hfill $\square$
\vspace{3mm}

Note that a by-product of Theorem \ref{TT1} is that for all
$ \varphi\in
C^\infty_0(\BbR^3, \BbC^2)$, and all $\,\nu\in [0, 1]$, the following Hardy-type
inhomogeneous inequality holds
$$\nu\int_{\BbR^3}\frac{|\varphi|^2}{|x|}\,
+\,\sqrt{1-\nu^2}\int_{\BbR^3}{|\varphi|^2}\,\leq\,
\int_{\BbR^3}\frac{|(\sigma\cdot\nabla)\varphi|^2}{\frac{\nu}{|x|}+1+\sqrt{1-\nu^2}}\,+\,
\int_{\BbR^3}{|\varphi|^2}\,.$$
This is just the inequality $Q^T_{\sqrt{1-\nu^2},\nu,0}(\varphi)\geq 0\,$
in the case $0\leq\nu<1$, and the case $\nu=1$ is obtained by passing to the limit.

\vspace{2mm}
Moreover, taking $\,\nu=1$ and functions $\,\varphi\,$ which concentrate near the
origin, the above inequality yields, in the limit, the following homogeneous one :
$$\int_{\BbR^3}\frac{|\varphi|^2}{|x|}\,dx\leq
\int_{\BbR^3}|x||{(\sigma\cdot\nabla)\varphi|^2}\,dx\quad\hbox{ for all }\quad 
\varphi\in
C^\infty_0(\BbR^3, \BbC^2)\,.$$

Actually, taking $\phi={\varphi\over \vert x\vert^{1/2}}$, this inequality is a direct
consequence of the standard Hardy inequality
$$\int_{\BbR^3} {\vert \phi\vert^2\over \vert x\vert^2}\,\leq\,
4\int_{\BbR^3} \vert \nabla\phi\vert^2\,=\,
4\int_{\BbR^3} \vert (\sigma\cdot\nabla)\phi\vert^2 \;.$$

\vspace{2mm}
\subsection{The min-max associated with the free-energy projectors.}

Here we define the splitting of ${\mathcal H}$ as follows: 
${\mathcal H}={\mathcal H}^f_+\oplus{\mathcal H}_-^f$, with \break
$\,{\mathcal H}_\pm^f=
\Lambda^f_\pm {\mathcal H}$, where
$$\Lambda^f_+=\chi_{(0, +\infty)}(H_0)=
\frac{1}{2}\left(\un+\frac{H_0}{\sqrt{1-\Delta}}\right)\,,$$
$$\Lambda^f_-=\chi_{(-\infty, 0)}(H_0)=
\frac{1}{2}\left(\un-\frac{H_0}{\sqrt{1-\Delta}}\right)\,.$$

As in Subsection 4.1, assumptions $(i)$ and $(ii)$ of Theorem \ref{S1} are satisfied,
with the same choice $\,a=c_2-1$.

\vspace{2mm}
With the new splitting ${\mathcal H}^f_\pm$, and the operator $\,A=H_0+V\,$, the min-max
values given by formula (\ref{min-max}) will be denoted by $\,\lambda_k^f(V)$.
This min-max principle based on free-energy projectors  was first introduced in
\cite{[Es-Se]}. Using some inequality proved in \cite{[Bur-Ev]} and
\cite{[Tix]},  we proved in \cite{[Dol-Est-Ser]} that  $\,\lambda_k^f(V)$
is indeed equal to the eigenvalue $\mu_k(V)$ for all potentials $V$
satisfying $\,-\frac{\nu}{|x|}\leq V\leq  0\,$, and all
$\,0\leq\nu<2\left(\frac{\pi}{2}+\frac{2}{\pi}\right)^{-1}\!\!\sim 0,9$. Here, we extend
this result to cover all $\,0\leq\nu<1$, and we obtain new inequalities as a by-product.

\vspace{2mm}
The main result of this subsection is the following

\begin{theorem}\label{TT67} Let $V$ a scalar potential satisfying
(\ref{V1})-(\ref{V2})-(\ref{V3}).
Then,  for all $\,k\geq 1$,  
\bq\label{free}\lambda^f_k( V)=\mu_k( V)\,.\eq
\end{theorem}

\noindent{\it Proof :} As in Subsection 4.1, we just have to consider the Coulomb
potential $\,V_\nu$, for $\,\nu\in [0, 1)$.

\vspace{2mm}
{\sl First Step :} Let $\varepsilon>0$ fixed and $\,V_{\nu, \varepsilon}\,$ as before.
Thanks to (\ref{import}),
Theorem \ref{TT5} applies to the one-parameter family $\,\nu\in [0,1)\to
A_{\nu, \varepsilon}:= H_0+V_{\nu, \varepsilon}\,$ with the projectors $\Lambda^f_\pm\,,$
and $a_-=-1\,,\;a_+=0\,$. So we get
$$\lambda^f_1(V_{\nu, \varepsilon})= \mu_1(V_{\nu, \varepsilon})\geq 0\,,$$
for all $\,\nu\in [0, 1)$. By Lemma \ref{N1}, this may be written
$$Q^f_{0,\nu, \varepsilon}(\psi_+)\geq 0\,,\quad\hbox{for all}\quad
\psi_+\in F^f_+:=\Lambda^f_+\left(C^\infty_0(\BbR^3,\BbC^{\,4})\right)\,,$$
with 
\begin{eqnarray}\label{3et}&\hspace{-4cm} Q^f_{E,\nu, \varepsilon}(\psi_+)=
||\psi_+||^2_{H^{1/2}}-(\psi_+, (E-V_{\nu, \varepsilon})\psi_+)\\&\hspace{11mm}+
\left(\Lambda^f_-|V_{\nu, \varepsilon}|\psi_+, \left(\Lambda^f_-(\sqrt{1-\Delta}+E+
|V_{\nu, \varepsilon}|)\Lambda^f_-\right)^{-1}\Lambda^f_-
|V_{\nu, \varepsilon}|\psi_+\right).\nonumber\end{eqnarray}

{\sl Second step :} Passing to the limit $\,\varepsilon\to 0\,$ in (\ref{3et}) with
$\psi_+$ and $\nu$ fixed, we get
\bq\label{4et}Q^f_{0,\nu, 0}(\psi_+)\geq 0\;,\quad\ \, \psi_+\in F^f_+
\eq
for all $\,\nu\in [0, 1)$. Then, applying Theorem \ref{S1} to $H_0+V_\nu$, we
obtain (\ref{free}), and the theorem is proved.
\hfill $\square$ 
\vspace{3mm}

Finally, note that some
inequalities can be derived from the free-energy min-max principle,
as in the Talman case:
for all $\,\nu\in [0, 1]$ and all
functions \break
$\,\psi_+\in \Lambda^f_+\left(C^\infty_0(\BbR^3,\BbC^{\,4})\right)$, we have

\begin{eqnarray*}
&\hspace{-55mm}\nu\!\Int_{\BbR^3}\!\!
\frac{|\psi_+|^2}{|x|}\,dx+\sqrt{1-\nu^2}\!\!\int_{\BbR^3} |\psi_+|^2\,dx
\\& \\&
\hspace{-70mm}\leq
\!\int_{\BbR^3} (\psi_+, \sqrt{1-\Delta}\,\psi_+)\,dx
\\
&  \\&+\,\nu^2\!\!\Int_{\BbR^3}\!\!
\left(\Lambda^f_-\left(\frac{\psi_+}{|x|}\right)\!, \!
\left(\Lambda^f_-\left(\!\!\sqrt{1-\Delta}\!\!+\tfrac{\nu}{|x|}\!+\!
\sqrt{1-\nu^2}
\right)\Lambda^f_-\right)^{-1}\!\!\!\!\!\!\Lambda^f_-\left(\tfrac{\psi_+}{|x|}
\right)\right)dx.
\end{eqnarray*}

Moreover, taking functions with support near the origin, we find, after rescaling
and passing to the limit,
a new homogeneous Hardy-type inequality. This inequality involves the projectors
associated with the  zero-mass free Dirac operator:

$$\Lambda^{f,0}_\pm:=\frac{1}{2}\left(\un\pm\frac{\alpha\cdot \hat{p}}{|
\hat{p}|}\right)\,,\;
 \hat{p}:=-i\nabla\,.$$
It may be written as follows :

\begin{eqnarray}&\hspace{-51mm} \Int_{\BbR^3} \frac{|\psi_+|^2}{|x|}\,dx\leq
\int_{\BbR^3} (\psi_+, | \hat{p}|\psi_+)\,dx\nonumber\\
&\hspace{5mm}+\Int_{\BbR^3}
\left(\Lambda^{f,0}_-\left(\frac{\psi_+}{|x|}\right), 
\left(\Lambda^{f,0}_-\left(| \hat{p}|+\frac{1}{|x|}
\right)\Lambda^{f,0}_-\right)^{-1}\Lambda^{f,0}_-\left(\frac{\psi_+}{|x|}\right)\right)\,dx\,,
\nonumber\end{eqnarray}
 for all
$\,\psi_+\in \Lambda^{f,0}_+\left(C^\infty_0(\BbR^3,\BbC^{\,4})\right).$

\vspace{2mm}
These two inequalities look like the ones obtained by Evans-Perry-Sieden\-top
\cite{[Eva-Per-Sie]}, Tix \cite{[Tix]} and Burenkov-Evans \cite{[Bur-Ev]},
but they are not the same. We do not know whether they can be obtained
by direct computations, as was the case in those works.

\newpage


\begin{center}\large\textsf{Erratum:\\ \emph{On the eigenvalues of operators with gaps.\\ Application to Dirac operators}}\end{center}

\small{{\sc Abstract}. In this \emph{erratum}, we address some closability issues that were ignored in~\cite{MR1761368}.}\medskip

\noindent {\bf AMS Subject Classification (2020):} Primary: 47B25. Secondary: 47A75, 49R50, 81Q10.
\vskip 0.5cm\noindent {\bf Keywords:} Variational methods; self-adjoint operators; quadratic forms; spectral gaps; eigenvalues; min-max principle; Rayleigh-Ritz quotients; Dirac operators.
\bigskip

In~\cite{MR4056270}, L.~Schimmer, J.P.~Solovej, and S.~Tokus construct a distinguished self-adjoint extension of a general symmetric operator with a gap and give a variational characterization of its eigenvalues, thus connecting the extension problem considered in~\cite{EstLos-07,EstLos-08} to the min-max principle for eigenvalues of self-adjoint operators studied in~\cite{MR1761368}. They also point out and solve several questions of closability and domain invariance that were not properly addressed in these papers. The min-max result of~\cite[Theorem~1.1]{MR4056270}, when applied to already self-adjoint operators, is analogous to the result of~\cite[Theorem~1.1]{MR1761368}, but with different assumptions. On the one hand, L.~Schimmer, J.P.~Solovej, and S.~Tokus do not assume that the subspace $F$, in which the min-max principle is defined, is a core. On the other hand, they strengthen an assumption denoted (i) in~\cite{MR1761368,MR4056270} and add a condition that we call $(\mathcal C)$ in this erratum. Moreover they suggest that $(\mathcal C)$ is also needed for~\cite[Theorem~1.1]{MR1761368}. They are completely right: the proof in~\cite{MR1761368} overlooks several closability issues, but can be completed under Condition~$(\mathcal C)$, and this is done in~\cite{MR4056270}. However, under minor corrections which are exposed below, the result of~\cite[Theorem~1.1]{MR1761368} also holds without assuming $(\mathcal C)$.\bigskip

\subsection*{Min-max characterization of the eigenvalues in a gap}\label{Sec:2001}

We first recall (up to a minor change that will be commented below) the assumptions of~\cite[Theorem~1.1]{MR1761368}.

Let $\mathcal H$ be a Hilbert space with norm $\|\cdot\|_{\scriptscriptstyle\mathcal H}$ and scalar product $(\cdot,\cdot)$ and let $A$ be a self-adjoint operator on $\mathcal H$ with domain $D(A)$. Let $\mathcal H_+$ and $\mathcal H_-$ be two orthogonal Hilbert subspaces of $\mathcal H$ such that $\mathcal H=\mathcal H_+\oplus\mathcal H_-$. Let $\Lambda_+$ and $\Lambda_-$ be the projectors on $\mathcal H_+$ and
~$\mathcal H_-$. We assume the existence of a core~$F$, \emph{i.e.}, a dense subspace of $D(A)$, such that :
\begin{itemize}
\item[(i)] $F_+=\Lambda_+F$ and $F_-=\Lambda_-F$ are two subspaces of $D(A)\,$,
\item[(ii)] $\displaystyle a=\sup_{x_-\in F_-\setminus\{0\}}\frac{(x_-,Ax_-)}{\|x_-\|^2_{\scriptscriptstyle\mathcal H}}<+\infty\,$,
\item[(iii)] $\lambda_1>a\,$.
\end{itemize}
Here $\lambda_1$ is the first of the min-max levels
\[
\lambda_k=\,\inf_{\scriptstyle\begin{array}{c}\mbox{\scriptsize $V$ subspace of $F_+$}\\[-4pt] \mbox{\scriptsize dim$\,V=k$}\end{array}}\,\sup_{ x\in(V\oplus F_-)\setminus\{0\}}\,\frac{(x,Ax)}{\|x\|^2_{\scriptscriptstyle\mathcal H}}\,,\quad k\ge1\,.
\]
Let $b=\inf\big(\sigma_{\rm ess}(A)\cap(a,+\infty)\big)\in[a,+\infty]$. For $k\ge 1$, $\mu_k$ denotes the $k^{\rm th}$ eigenvalue of~$A$ in the interval $(a,b)$, counted with multiplicity, \emph{if this eigenvalue exists}. If there is no $k^{\rm th}$ eigenvalue, we take $\mu_k = b\,$.
\begin{theorem}\label{E-S1} With the above notations, and under assumptions {\rm(i)-(ii)-(iii)}, we have
\[
\lambda_k\,=\,\mu_k\,,\quad\forall\,k\ge1
\]
and, as a consequence, $b\,=\lim_{k\to+\infty}\lambda_k\,=\,\sup_{k\ge1}\lambda_k\,>\,a\,$.
\end{theorem}

Compared to~\cite[Theorem~1.1]{MR1761368}, the only change is that in (i) one has to assume that $F_+$ and $F_-$ are subspaces of $D(A)$ and not of the form domain $\mathcal F(A)$. This is weaker than Condition (${\mathcal C}$)
\begin{center}
\emph{The operator $\Lambda_-A_{|F_-}:F_- \to\mathcal H_-$ is essentially self-adjoint,}
\end{center}
which is called assumption (iii) in~\cite[Theorem~1.1]{MR4056270}. Note that a similar change is also needed in the assumptions of the continuation principle of~\cite[Theorem~3.1]{MR1761368}: in hypothesis (j), $\mathcal{F}(A_0)$ has to be replaced by $\mathcal{D}(A_0)$.

Note that \cite[Theorem~1.1]{MR1761368} supposedly applies to Dirac-Coulomb operators $A={\bf \alpha}\cdot {\bf p}+{\beta}\,m -\frac{\nu}{r}$ with $F=C^\infty_c(\mathbb{R}^3,\mathbb{C}^4)$ for any $0\leq \nu<1$. Such a claim is incorrect, since $F$ is not a core of $A$ when $\sqrt3/2<\nu<1$. This issue was pointed out to us by several researchers, and a partial solution, that requires the replacement of $C^\infty_c$ by $H^{1/2}$, was given by S.~Morozov and D.~M\"{u}ller in~\cite{MorMul-15,Mueller-2016}. However, the issue can also be solved in $C^\infty_c$ by first applying the theorem to regularized Dirac-Coulomb operators $A_\varepsilon={\bf \alpha}\cdot {\bf p}+{\beta}\,m -\frac{\nu}{r+\varepsilon}$, and then passing to the limit in the norm-resolvent sense as $\varepsilon\to 0$, as done in~\cite{MR3960263}. Of course, another possibility is to apply directly the result of~\cite{MR4056270}, which does not make use of the assumption that $F$ is a core.

\subsection*{Changes in the proof of Theorem 1.1 of~\cite{MR1761368}}\label{Sec:Changes}

As pointed out in~\cite{MR4056270}, several closability issues were overlooked in the proof of~\cite[Theorem~1.1]{MR1761368}. Let us list the changes that have to be done in the proof to make it correct. The numbering of the formulae refers to~\cite{MR1761368}.

In~\cite{MR1761368} we introduced the norm
\[
N (y_-) = \sqrt{(a+1)\|y_-\|^2_{_{\mathcal H}} - (y_-, Ay_-)}
\]
on $F_-$ and claimed without proof that the completion of $F_-$ for this norm can be identified with a subspace of ${\mathcal H}_-$. In other words, we claimed that the quadratic form $q_-(y_-)=-(y_-, Ay_-)$ is closable in ${\mathcal H}_-$. Unfortunately we cannot prove this closability under the assumption $F_-\subset {\mathcal F}(A)$, but this becomes a standard fact if one assumes that $F_-\subset {\mathcal D}(A)$ as in~\cite{MR4056270}. Indeed, the operator $-\Lambda_-A_{|F_-}$ is then a symmetric and bounded from below operator defined on ${\mathcal H}_-$ with domain $F_-$. One can use the Friedrichs extension theorem~\cite[Theorem~X.23]{MR0493420} to establish that the form $q_-$ is closable in ${\mathcal H}_-$. Denoting by~$\overline{F}_-^N$ the form-domain of its closure $\overline{q}_-$, there is a unique self-adjoint operator $B$ with domain ${\mathcal D}(B)\subset \overline{F}_-^N$ such that for any $y_-\in {\mathcal D}(B)$, $\overline{q}_-(y_-)=(y_-,By_-)$. The operator $B$ is the Friedrichs extension of $-\Lambda_-A_{|F_-}$.

The operator $L_E: F_+\to {\mathcal H}_-$ is defined by
\[
L_E x_+=(B+E)^{-1}\Lambda_-Ax_+\,.
\]
In order to define $L_E$ without the assumption $F_-\subset {\mathcal F}(A)$, a delicate construction was needed in~\cite{MR1761368}, while the definition of $L_E$ is now straightforward.

Then we considered the completion $X$ of $F_+$ for the norm $n_E(x_+)=\|x_+ + L_E x_+ \|_{_{\mathcal H}}$. We claimed without proof that $X$ can be identified with a subspace of ${\mathcal H}_+$. In other words, we claimed that the operator $L_E$ is closable in ${\mathcal H}_+$. As pointed out in~\cite{MR4056270}, this is far from obvious, and probably wrong without an additional assumption. In~\cite{MR4056270}, the authors prove that $L_E$ is closable assuming that Condition $(\mathcal C)$ holds. The main goal of the present \emph{erratum} is to explain that this additional assumption is not needed in the proof of Theorem~\ref{E-S1}.
 
Without condition $(\mathcal C)$, we cannot claim that $L_E$ is closable, but instead of $X$, we can consider the closure $\overline{\Gamma}_E$ of its graph represented by $\Gamma_E=\big\{x_+ + L_E x_+:\, x_+ \in F_+\big\}$ in the Hilbert space ${\mathcal H}$ identified with ${\mathcal H}_+\times {\mathcal H}_-$. The closed subspace $\overline{\Gamma}_E$ does not necessarily represent a graph over ${\mathcal H}_+$, but this does not affect too much the remainder of the proof. We just need to modify the definitions of some mathematical objects. Several expressions that were defined as functions of $x_+\in F_+$ are now considered as functions of $x_+ +L_Ex_+\in \Gamma_E$ and their extensions by density become functions of $x\in \overline{\Gamma_E}$. In particular, $n_E$ (resp. $\overline{n}_E$) becomes the restriction of $\Vert\cdot\Vert_{\mathcal{H}}$ to $\Gamma_E$ (resp. $\overline{\Gamma_E}$); in formula~(7) and everywhere else in the sequel, $Q_E$ becomes a map from $\Gamma_E$ to $\mathbb R$ while $Q_E(x_+)$ has to be replaced by
\[
Q_E\big(x_++L_Ex_+\big)=\big((A-E)x_+,x_+) + (L_E x_+,(B+E) L_E x_+\big)\,;
\]
in formula~(8), $Q_E(\Lambda_+x)$ becomes $Q_E(\Lambda_+x+L_E\Lambda_+x)\,$. Formulae (10) and (11) have to be rewritten as
\be{10'}\tag{10'}
\Vert x_+\Vert_{\mathcal H}\leq\Vert x_++L_{E'}x_+\Vert_{\mathcal H}\leq\Vert x_++L_Ex_+\Vert_{\mathcal H}\leq\frac{E'-a}{E-a}\,\Vert x_++L_{E'}x_+\Vert_{\mathcal H}
\ee
and
\begin{multline}\label{11'}\tag{11'}
(E'-E)\,\Vert x_++L_{E'}x_+\Vert_{\mathcal H}^2\leq Q_E\big(x_++L_Ex_+\big)-Q_{E'}\big(x_++L_{E'}x_+\big)\\
\leq (E'-E)\,\Vert x_++L_Ex_+\Vert_{\mathcal H}^2\,.
\end{multline}
Up to these changes,~\cite[Lemma~2.1]{MR1761368} remains as stated previously. Note that formula~\eqref{10'} implies that the map $x_++L_Ex_+\mapsto x_++L_{E'}x_+$ is an isomorphism between $\Gamma_E$ and $\Gamma_{E'}$ which extends to an isomorphism between $\overline \Gamma_E$ and $\overline \Gamma_{E'}$. The replacement of $x_+\in F_+$ by $x_+ +L_Ex_+\in \Gamma_E$ has also to be done in the definition of the norm ${\mathcal N}_E$, which becomes
\[
{\mathcal N}_E(x_++L_Ex_+)=\sqrt{Q_E(x_++L_Ex_+)+(K_E+1)\Vert x_++L_Ex_+\Vert_{\mathcal H}^2}\,.
\]
As a consequence, formula~(13) takes a slighly different form, since ${\mathcal N}_E$ and ${\mathcal N}_{E'}$ are now defined on different spaces. It becomes
\begin{multline}\label{13'}\tag{13'}
c(E,E')\,{\mathcal N}_E(x_++L_Ex_+)\leq {\mathcal N}_{E'}(x_++L_{E'}x_+)\leq C(E,E')\,{\mathcal N}_E(x_++L_Ex_+)\,,\\
\forall x_+ \in F_+\,.
\end{multline}
Let $\Pi_E:\,\mathcal{H}\to \overline{\Gamma}_E$ be the orthogonal projector on $\overline{\Gamma}_E$. Note that for any $x\in\Gamma_E$ we have $Q_E(x)=(x, S_E x)$ where
\[
S_E x=\Pi_E\big(\Lambda_+(A-E) \Lambda_+x +(B+E)\Lambda_-x\big)\,.
\]
It is clear that $S_E$ is a symmetric operator bounded from below on the Hilbert space $\overline{\Gamma_E}$, with domain $\Gamma_E$. So~$Q_E$ is closable in $\overline{\Gamma_E}$. We denote its closure $\overline{Q}_E$. Its domain is denoted by $G_E$ and the corresponding extended norm is denoted by $\overline{\mathcal{N}}_E$. Note that $G_E$ is a subspace of $\overline{\Gamma_E}$, so it depends on $E$, while in~\cite{MR1761368} its analogue $G$ was considered as a subspace of $\mathcal{H}_+$, independent of $E$ and, as a consequence, the norms $\overline{\mathcal{N}}_E$ were all equivalent. Now~\eqref{13'} implies that any two normed spaces $(G_E,\overline{\mathcal{N}}_E)$ and $(G_{E'},\overline{\mathcal{N}}_{E'})$ are isomorphic.

We denote by $T_E$ the Friedrichs extension of $S_E$. It is a self-adjoint operator in $\overline{\Gamma_E}$ with domain $\mathcal{D}(T_E)$ and we have $\Gamma_E\subset \mathcal{D}(T_E)\subset G_E\subset \overline{\Gamma_E}\subset\mathcal{H}$. In~\cite{MR1761368}, $\mathcal{D}(T_E)$ was a subspace of~$\mathcal{H}_+$, but this is no longer true with our new definition. We do not know whether $\mathcal{D}(T_E)$ is a graph over $\mathcal{H}_+$, but the important point is that the graph $\Gamma_E$ is a form-core of $T_E$ and~$G_E$ is its form-domain. In formula~(14) for the min-max levels, $G$ should be replaced by $G_E$ and the notation~$x_+$ should be replaced by $x$, since this variable no longer belongs to~$\mathcal{H}_+$. Moreover, for any $x\in\overline \Gamma_E$, the extended norm $\overline n_E(x)$ considered in~\cite{MR1761368} is replaced by $\Vert x\Vert_{\mathcal{H}}$. Formula (14) thus becomes
\be{14'}\tag{14'}
\ell_k(T_E)=\,\inf_{\scriptstyle\begin{array}{c}\mbox{\scriptsize $V$ subspace of $G_E$}\\[-4pt] \mbox{\scriptsize dim$\,V=k$}\end{array}}\,\sup_{x\in V\setminus\{0\}}\;\frac{\overline Q_E(x)}{\|x\|_{\mathcal H}^2}\,.
\ee
In formula~(15) of~\cite[Lemma~2.2]{MR1761368}, $x_+$ should be replaced by $x_++L_\lambda x_+$. The lemma remains otherwise unchanged. In part (a) of its proof, the notation $Q_\lambda(x_+)$ is used repeatedly. It should be replaced everywhere by $Q_\lambda(x_++L_\lambda x_+)$, but no other change has to be done.

Similarly, in the sequel of the proof of~\cite[Theorem~1.1]{MR1761368}, $T_{\lambda_k}x_+$ should be replaced by $T_{\lambda_k}(x_++L_{\lambda_k}x_+)$ and $G'$ by $G_{\lambda_k}'$; 
$\mathcal{A}$ is now the polar form of the quadratic form $\Gamma_{\lambda_k}\ni\tilde{y} \mapsto Q_{\lambda_k}(\tilde{y})+\lambda_k\,\Vert\tilde{y}\Vert^2_{\mathcal{H}}\,$; one should replace $Q_{\lambda_k}(\Lambda_+ y)$ by $Q_{\lambda_k}(\tilde{y})$ in formula~(20) and its proof.

Moreover we recall that $\Gamma_E$ is a form-core of~$T_E$ and $G_E$ is its form-domain, so that~\eqref{14'} is equivalent to
\[
\ell_k(T_E)=\,\inf_{\scriptstyle\begin{array}{c}\mbox{\scriptsize $V_+$ subspace of $F_+$}\\[-4pt] \mbox{\scriptsize dim$\,V_+=k$}\end{array}}\,\sup_{x_+\in V_+\setminus\{0\}}\;\frac{ Q_E(x_++L_Ex_+)}{\|x_++L_Ex_+\|_{\mathcal H}^2}\,.
\]

\subsection*{An example}\label{Sec:Example}

The additional assumption $F_\pm\subset \mathcal{D}(A)\,$ seems harmless: in all examples of the literature we are aware of, it is satisfied. Condition ($\mathcal{C}$) is satisfied in many situations of interest in physics, as explained in~\cite{MR4056270}. But in this short section we give an example where ($\mathcal{C}$) is not satisfied. Let $V$ be an electric potential in $\mathbb R^3$ of the form
\[
V(x)=-\frac{\nu_1}{\vert x \vert}+\frac{\nu_2}{\vert x-x_0 \vert}
\]
with $0<\nu_1<\sqrt3/2$, $3/4<\nu_2<\sqrt3/2$ and $x_0\in\R^3\setminus\{0\}$. On the Hilbert space $\mathcal{H}=\mathrm L^2(\mathbb{R}^3,\mathbb{C}^4)$, consider the associated Dirac-Coulomb operator
\[
A={\bf \alpha}\cdot {\bf p}+{\beta}\,m +V
\]
with domain $\mathrm H^1(\mathbb{R}^3,\mathbb{C}^4)$. Since $0\leq \nu_1\,,\,\nu_2<\sqrt3/2\,,$ $A$ is self-adjoint and admits the subspace $F=C^\infty_c(\mathbb{R}^3,\mathbb{C}^4)$ as a core. Let $\Lambda_\pm={\bf 1}_{\mathbb{R}_\pm}({\bf \alpha}\cdot {\bf p}+{\beta}\,m)$. As was proved by C.~Tix in~\cite{MR1608118}, the Brown-Ravenhall operators $\Lambda_+ A_{|\Lambda_+F}$ and $-\Lambda_- A_{|\Lambda_-F}$ are both positive, so that one can apply Theorem~\ref{E-S1} with $a=0$, but as shown in \cite[Corollary~3]{tix1997self}, $\Lambda_- A_{|\Lambda_-F}$ is not essentially self-adjoint, since $\nu_2 >3/4$.

\medskip\noindent{\bf Acknowledgment:} We thank the authors of~\cite{MR4056270} for letting us know about the closability issue in~\cite{MR1761368}.

\end{article}
\end{document}